\documentclass[11pt,reqno]{amsart}
\usepackage {amsmath,amsfonts,enumerate,  epsfig,   indentfirst,  graphicx,euscript, amssymb,amsthm,amscd}

\usepackage[dvips]{psfrag}

\usepackage[latin1]{inputenc}

\makeindex

\usepackage{amssymb}
\usepackage[latin1]{inputenc}
\usepackage{enumerate}
\usepackage{indentfirst}
\usepackage{euscript}
\usepackage{graphicx}
\usepackage{amsmath}
\usepackage{amsfonts}
\usepackage{amssymb}
\makeindex

\makeindex
\theoremstyle{plain}
\newtheorem{proposition}{Proposition}

\newtheorem{corollary}{Corollary}
\newtheorem{remark}{Remark}
\theoremstyle{definition}

\newtheorem{theorem}{Theorem}

\newtheorem{problem}[theorem]{Problem}

\title[]
{Surfaces Around Closed Principal Curvature Lines, an Inverse
Problem}

\author{R. Garcia, L. F. Mello and J. Sotomayor}

\begin{document}
\maketitle

\begin{abstract}
Given a non circular spacial closed curve whose total torsion is an integer multiple of
$2\pi$, we construct a germ of a smooth surface that contains it as a hyperbolic principal
cycle.
\end{abstract}

\section{Introduction}

Let $\alpha:\mathbb M \to \mathbb R^3$ be a $ C^{r}$ immersion of a
smooth, compact and oriented, two--dimensional manifold $\mathbb M$
into space ${\mathbb R}^{3}$ endowed with the canonical inner
product $<.,.>$. It  will be assumed that $r\geq 4$.

The {\em Fundamental Forms} of $\alpha $ at a point $p$ of ${\mathbb
M}$ are the symmetric bilinear forms on ${\mathbb T}_p \mathbb M$
defined as follows \cite{spi}, \cite{st}:
\[
\aligned I_{\alpha }(p;v,w)=& \left< D\alpha
(p;v),D\alpha (p;w) \right>,\\
II_{\alpha}(p;v,w)=& \left< -DN_{\alpha}(p;v),D\alpha(p;w)
\right>.\endaligned
\]
Here, $N_{\alpha }$ is the positive normal of the immersion
$\alpha$.

The first fundamental form in a local chart $(u,v)$ is defined by
$I_\alpha= Edu^2+2Fdudv+Gdv^2$, where $E=\left<\alpha_u,
\alpha_u\right>$, $F=\left<\alpha_u, \alpha_v\right>$ and
$G=\left<\alpha_v, \alpha_v\right>$.

The second fundamental form relative to the unitary normal vector
$N_\alpha= ( \alpha_u\wedge \alpha_v)/|\alpha_u\wedge \alpha_v|$ is
given by $II_\alpha= edu^2+2fdudv+gdv^2$, where

$$ \text{\small $ e=\frac{  \det[ \alpha_u,\alpha_v,\alpha_{uu}]}{
\sqrt{   EG -F^2}}, \;
 f=\frac{  \det[ \alpha_u,\alpha_v,\alpha_{uv}]}{\sqrt{EG-F^2}},\;
g=\frac{ \det[ \alpha_u,\alpha_v,\alpha_{vv}]}{\sqrt{ EG-F^2 }}.$}  $$

In a local chart $(u,v)$ the principal directions of an immersion
$\alpha$ are defined by the implicit differential equation
\begin{equation}\label{eq:p}
(Fg-Gf)dv^2+(Eg-Ge)dudv+(Ef-Fe)du^2=0.
\end{equation}

The {\em umbilic set} of $\alpha$, denoted by ${\mathcal U}_\alpha$,
consists on the points where the three coefficients of equation
(\ref{eq:p}) vanish simultaneously.

The regular integral curves of equation \eqref{eq:p} are called {\em
principal curvature lines}. This means curves $ c(t)=(u(t), v(t))$,
differentiable on an interval, say $J$, with non--vanishing tangent
vector there, such that, for every $t \in J$, it holds that
\[
\aligned (Fg-Gf)(u(t),v(t)) \left( \frac{dv(t)}{dt} \right)^2+(Eg-Ge)(u(t),v(t))\frac{du(t)}{dt}\frac{dv(t)}{dt}\\
+(Ef-Fe)(u(t),v(t)) \left( \frac{du(t)}{dt} \right)^2  = 0
\endaligned
\]
and $c(J)\cap{\mathcal U}_\alpha=\emptyset$.

When the surface $\mathbb M$ is oriented, the principal curvature
lines on $\mathbb M \setminus {\mathcal U}_\alpha$  can be assembled
in two one--dimensional orthogonal foliations which will be denoted
by ${\mathcal F}_1(\alpha)$ and ${\mathcal F}_2(\alpha)$. Along the
first (resp. second), the normal curvature $II_{\alpha}(p)$  attains
its minimum $k_1 (p)$, denominated the {\em minimal principal
curvature at} $p$, (resp. maximum $k_2 (p)$, denominated the {\em
maximal principal curvature at} $p$).

The triple ${\mathcal P}_\alpha=\{ {\mathcal F}_1(\alpha),{\mathcal
F}_2(\alpha), {\mathcal U}_\alpha\}$ is called the {\it principal
configuration} of the immersion $\alpha$, \cite{gs1}, \cite{gs2}. For a
survey about the qualitative theory of principal curvature lines see
\cite{arxiv}.

A closed principal curvature line,
also called a {\em principal cycle}, is called {\em  hyperbolic} if
the first derivative of the Poincar\'e return map associated to it
is different from one.

In \cite{li} and \cite{pc} it was proved that a regular closed line of
curvature on a surface has as total torsion a multiple of $2\pi$.
In this paper we consider the following inverse problem.

\begin{problem} \label{p:1}
Given a simple closed Frenet curve, that is a smooth regular curve of $\mathbb R^ 3$ with non zero curvature,
is there an oriented embedded surface that contains it as a hyperbolic principal cycle?
\end{problem}

It will be shown that this Problem  has a positive answer in the
case that the curve is a Frenet, non circular, curve such that its total
torsion is an integer multiple of $2\pi$.

The interest of hyperbolic principal cycles is that the asymptotic
behavior of the principal foliation around them is determined. The
first examples of hyperbolic principal cycles on surfaces were
considered by Gutierrez and Sotomayor in \cite{gs1}, where their
genericity and structural stability  were also established.

\section{Preliminary Results}

Let $c:[0,L] \to \mathbb R^3$ be a smooth {\em simple, closed,
regular curve} in $\mathbb R^3$ with positive curvature $k$ and of
length $L>0$, i.e., a Frenet curve. Let also $\mathbf{c}=c([0,L])$. Consider the Frenet
frame $\{t,n,b\}$ along $\mathbf{c}$ satisfying the equations
\begin{equation}\label{eq:fre}\aligned
t^\prime(s)=& k(s) n(s),\\
n^\prime(s)=&-k(s)t(s)+\tau(s) b(s),\\
b^\prime(s)= &- \tau(s) n(s).\endaligned
\end{equation}
Here $k>0$
is the curvature and $\tau$ is the torsion of $\mathbf{c}$.

Consider the parametrized surface of class $C^r, \; r\geq 4$, defined by the equation
\begin{equation}\label{eq:a}\aligned
\alpha(s,v)=&  c(s)+\left[ \cos \theta(s) n(s)+ \sin\theta(s) b(s) \right] v \\
+& \left[ \cos \theta(s) b(s)- \sin\theta(s) n(s) \right] \left[
\frac 12 A(s) v^2+ \frac 16 B(s)v^3+v^4C(s,v) \right]\\
= &c(s)+v (N\wedge T)(s)+ \left[ \frac 12 A(s) v^2+ \frac 16
B(s)v^3+v^4C(s,v) \right] N(s).\endaligned
\end{equation}

\begin{figure}[htbp]
\begin{center}
\psfrag{n}{$N$} \psfrag{t}{$T$} \psfrag{nt}{$N\wedge T$}
\psfrag{c}{$\mathbf{c}$} \psfrag{o1}{$0$} \psfrag{s}{$s$}
\psfrag{l}{$L$} \psfrag{v}{$v$} \psfrag{o}{$O=(0,0,0)$}
\psfrag{a}{$\alpha(s,v)$} \psfrag{a1}{$\alpha$}
\includegraphics[scale=0.5]{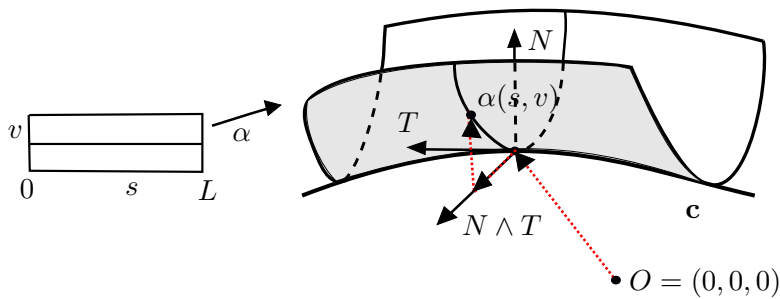}
\caption{Germ of a parametrized surface $\alpha(s,v)$ near a curve
$\mathbf{c}$. \label{fig:1}}
\end{center}
\end{figure}

For an illustration see Fig. \ref{fig:1} and \cite{coloquio}, \cite{gs1}.

Here, $\mathbf{c}^\prime(s)=t(s)=T(s)$, $\theta(s)=\theta(s+L)$, $A(s)=A(s+L)$, $B(s)=B(s+L)$, $C(s,v)
= C(s+L,v)$, $C(s,0)=0$, are smooth $L$--periodic functions with respect to $s$
and $v$ is small.

\begin{proposition} \label{prop:lc}
The curve $\mathbf{c}$ is the union of principal curvature lines of $\alpha$
if and only if
\begin{equation} \label{eq:tt}
\tau(s)+\theta^\prime(s)=0, \; \theta(0)=\theta_0, \;\;\; \int_0^L \tau(s)ds=2 m \pi, \;\; m\in \mathbb Z.
\end{equation}

Moreover, for any solution $\theta(s)$ of equation \eqref{eq:tt} the
parametric surface defined by equation \eqref{eq:a} is a regular,
oriented and embedded surface in a neighborhood of $\mathbf{c}$. The umbilic
set ${\mathcal U}_\alpha \cap \mathbf{c} $ is defined by the
equation
\begin{equation}\label{eq:umb}
A(s)+k(s)\sin\theta(s)=0.
\end{equation}
\end{proposition}

\begin{proof}
We have that $N(s)=N_\alpha(s,0)= \cos \theta(s) b(s)- \sin\theta(s)
n(s)$. By Rodrigues formula it follows that $\mathbf{c}$ is a principal
curvature line  (union of maximal and minimal principal lines) if
and only if $N^\prime(s) +\lambda(s) t(s)=0$. Here $\lambda$ is a
principal curvature (maximal or minimal).

Differentiating $N$ leads to
\[
N^\prime(s)=\sin\theta(s) k(s)
t(s)-[\tau(s)+\theta^\prime(s)][\cos\theta(s) n(s)+\sin\theta(s)
b(s)].
\]
Therefore $\mathbf{c}$ is a principal line
(union of maximal and minimal
principal lines
and umbilic points) if and only if $\tau(s)+\theta^\prime(s)=0$.

By the definition of $\alpha$ it follows that $\alpha_s(s,0)=t(s)$
and $\alpha_v(s,0)=\cos\theta(s) n(s)+\sin\theta(s) b(s) $ are
linearly independent and so by the local form of immersions it
follows that $\alpha$ is locally a regular surface in a
neighborhood of $\mathbf{c}$.

Since the total torsion is an integer multiple
of $2\pi$ and $\tau(s)+\theta^\prime(s)=0$ it follows that for any
initial condition $\theta(0)=\theta_0$ equation \eqref{eq:a}
defines an oriented and embedded surface containing $\mathbf{c}$ and having
it as the union of principal lines
and umbilic points.

Supposing that $\tau(s)+\theta^\prime(s)=0$  it follows that the
coefficients of the first and second fundamental forms of $\alpha$
are given by
\begin{equation}\label{eq:12}\aligned
E(s,v)=& 1-2k(s)\cos\theta(s) v \\
+& \left[ \frac 12 k(s)^2(1+\cos 2\theta(s))+ k(s) A(s)\sin\theta(s) \right] v^2 +O (v^3), \\
F(s,v)=& \frac 12A^\prime (s) A(s) v^3+ O(v^4),\\
G(s,v)= &1+A(s)^2 v^2+O(v^3),\\
e(s,v)= &-k(s)\sin\theta(s)+ k(s)\cos\theta(s)(\sin  \theta(s)-  A(s) ) v+ O(v^2),\\
f(s,v)=& A^\prime(s) v+O(v^2),\\
g(s,v)=& A(s)+B(s) v+ O(v^2).
\endaligned
\end{equation}

By equations \eqref{eq:p} and \eqref{eq:12} it follows that the
coefficients of the differential equation of principal curvature
lines are given by
\begin{equation}\label{eq:lmn}\aligned
L(s,v)=& (Fg-Gf)(s,v) = -A^\prime(s) v+ O(v^2),\\
M(s,v)=& (Eg-Ge)(s,v)= A(s)+k(s)\sin \theta(s) \\
+& \left[ B(s)-k(s)A(s) \cos\theta(s)-\frac 12 k(s)^2\sin 2\theta(s) \right] v+O(v^2),\\
N(s,v)=& (Ef-Fe)(s,v)= A^\prime(s) v+ O(v^2).
\endaligned
\end{equation}

The umbilic points along $\mathbf{c}$ are given by the solutions of
$M(s,0)=A(s)+k(s)\sin\theta(s) =0$ which corresponds to the equality
between the principal curvatures $k_1(s) =-k(s)\sin\theta(s)$ and
$k_2(s)=A(s)$.
\end{proof}

\begin{remark}\label{rm:3}
The one parameter family of surfaces
$M(\theta_0)=\alpha_{\theta_0}([0,L]\times
(-\epsilon,\epsilon))\setminus \mathbf{c}$ defined by equations
\eqref{eq:a} and \eqref{eq:tt} is a foliation of a neighborhood of
$\mathbf{c}$ after  $\mathbf{c}$ is removed.
For all $\theta_0$ the curve $\mathbf{c}$ is a principal cycle of
$\alpha_{\theta_0}$. This follows from the theorem of existence and
uniqueness of ordinary differential equations and smooth dependence
with initial conditions of $\alpha_{\theta_0}$ and boundary
conditions given by equation \eqref{eq:tt}.
\end{remark}

\section{Hyperbolic Principal Cycles}\label{sc:p}

In this section it will be given a solution to the problem
formulated in the Introduction.

Let $\alpha_{\theta_0}$ be the surface defined by equation
\eqref{eq:a} and associated to the Cauchy problem given by equation \eqref{eq:tt}.

\begin{theorem} \label{th:lp}
Consider the oriented parametric surface $\alpha_{\theta_0}$ of class $C^r, \;r\geq 4$, defined by
equations \eqref{eq:a} and \eqref{eq:tt} such that ${\mathcal
U}_{\alpha_{\theta_0}} \cap \mathbf{c}=\emptyset$. Then $\mathbf{c}$ is a
hyperbolic principal cycle of $\alpha_{\theta_0}$ if and only if
\begin{equation}\label{eq:ch}
\Lambda (\theta_0)=\int_0^L
\frac{A^\prime(s)}{A(s)+k(s)\sin\theta(s)}ds\ne 0.
\end{equation}
The coefficient $\Lambda (\theta_0)$ is called the characteristic
exponent of the Poincar\'e map associated to $\mathbf{c}$.
\end{theorem}

\begin{proof} The principal curvatures are given by
\[
\aligned k_1(s,v)=&-k(s)\sin\theta(s)+(k(s)\sin\theta(s)+A(s)) k(s)\cos\theta(s)v+O(v^2),\\
k_2(s,v)=& A(s)+B (s)v+O(v^2).\endaligned
\]

The first return map $\pi:\{s=0\}\to \{s=L\}$ defined by
$\pi(v_0)=v(L,v_0)$, with $v(0,v_0)=v_0$, satisfies the variational
equation
\[
M(s,0) v_{sv_0}(s)+ N_v(s,0)v_{v_0}(s)=0.
\]
By equation \eqref{eq:lmn} it  follows that
\[
-\frac{N_v}{M}(s,0)=  -\frac{A^\prime(s)}{A(s)+k(s)\sin\theta(s)}.
\]
Integration of the above equation leads to the result.
\end{proof}

\begin{remark}\label{r:1}
The criterium of hyperbolicity of a principal cycle was established
by Gutierrez and Sotomayor in \cite{gs1}, \cite{gs2}. They proved that a
principal cycle $\mathbf{c}$ is hyperbolic if and only if
\[
\int_\mathbf{c}\frac{dk_1}{k_2-k_1} =\int_\mathbf{c}
\frac{dk_2}{k_2-k_1}= \frac{1}2 \int_\mathbf{c} \frac{d{\mathcal
H}}{\sqrt{{\mathcal H}^2-{\mathcal K}}}\ne 0.
\]
Here ${\mathcal H}= (k_1+k_2)/2$ and ${\mathcal K}=k_1 k_2$ are
respectively the {\em arithmetic mean} and the {\em Gauss}
curvatures of the surface.
\end{remark}

\begin{proposition} \label{prop:d2} Consider the family of oriented parametric
surfaces $\alpha_{\theta_0}$ defined by equations \eqref{eq:a} and
\eqref{eq:tt} such that ${\mathcal U}_{\alpha_{\theta_0}} \cap
\mathbf{c}=\emptyset$ for all $\theta_0$. Then the following holds
\begin{equation}\label{eq:d2}
\Lambda^\prime(\theta_0)=
\int_0^L \frac{k(s) A^\prime (s)\cos(\theta_0-\int_0^s
 \tau(s)ds)}{[k(s)\sin(\theta_0-\int_0^s \tau(s)ds)+A(s)]^2}ds. \end{equation}

\end{proposition}

\begin{proof} Direct differentation of equation \eqref{eq:ch}.
\end{proof}

\begin{theorem} \label{th:lp2}
Let $\mathbf{c}$ be a smooth
curve, that is a closed Frenet curve of length
$L$ in $\mathbb R^3$ such that $\tau$ is not identically zero and
$\int_0^L \tau(s)ds=2m\pi, \; m\in \mathbb Z$. Then there exists a
germ of an oriented surface of class $C^r$, $r\geq 4$, containing $\mathbf{c}$ and having it as
a hyperbolic principal cycle.
\end{theorem}

\begin{proof} Consider the parametric surface defined by equation \eqref{eq:a}.
By Proposition \ref{prop:lc}, $\mathbf{c}$ is a principal cycle when
$\theta^\prime(s)= -\tau(s),\; \theta(0)=\theta_0,$ $ \int_0^L
\tau(s)ds=2m\pi$. Taking $A(s)=(1- \sin\theta(s))k(s)$ it follows
that $M(s,0)=k(s)>0$ and so ${\mathcal U}_\alpha \cap
\mathbf{c}=\emptyset$. So $\mathbf{c}$ is a closed principal line of the
parametric surface $\alpha$.

By Theorem \ref{th:lp} it follows that $\mathbf{c}$ is hyperbolic if and only
if
\[
\Lambda=\ln(\pi^\prime(0))=\int_0^L
\frac{A^\prime(s)}{A(s)+k(s)\sin\theta(s)}ds= - \int_0^L
\frac{[k(s)\sin\theta(s)]^\prime}{ k(s) } ds\ne 0.
\]

By assumption the function $ k(s)\sin\theta(s) $ is not constant. In
fact, $k(s)>0 $ and as $\tau$ is not identically
equal to zero it follows
that $\sin\theta(s)= \sin(\theta_0-\int_0^s \tau(s)ds)$
is not constant.

If $\Lambda\ne 0$ it follows that $\mathbf{c}$ is hyperbolic and
this ends the proof.

In the case when $k(s)\sin\theta(s)$ is not constant and $\int_0^L
\frac{[k(s)\sin\theta(s)]^\prime}{k(s)} ds=0$, consider the
deformation of $\alpha$ given by
\[
\alpha_\epsilon(s,v)=\alpha(s,v)+\epsilon \frac{a(s)}2 v^2 [\cos
\theta(s) b(s)- \sin\theta(s) n(s)], \;
a(s)=[k(s)\sin\theta(s)]^\prime.
\]
Then $\mathbf{c}$ is a principal cycle of $\alpha_\epsilon$ and the principal
curvatures are given by
\[
\aligned
k_1(s,0,\epsilon)= &-k(s)\sin\theta(s),\\
 k_2(s,0,\epsilon)=&A(s) +\epsilon  [k(s)\sin\theta(s)]^\prime
 =  (1- \sin\theta(s))k(s) +\epsilon  [k(s)\sin\theta(s)]^\prime. \endaligned
\]
Therefore by Theorem \ref{th:lp} and Remark \ref{r:1} it follows
that
{\small
\[
\Lambda(\epsilon)= \ln(\pi^\prime_{\epsilon}(0))=-\int_0^L\frac{
k_1^\prime(s,0,\epsilon)}{k_2(s,0,\epsilon)-k_1(s,0,\epsilon)}ds=
\int_0^L \frac{[k(s)\sin\theta(s)]^\prime}{ k(s) +\epsilon
[k(s)\sin\theta(s)]^\prime} ds.
\]}
Differentiating the above equation with respect to $\epsilon$ and
evaluating in $\epsilon=0$ it follows that
\[
\Lambda^\prime(0)=\frac{d}{d\epsilon}\big(\ln(\pi^\prime_{\epsilon}(0))\big)|_{\epsilon=0}
= \int_0^L \left[ \frac{  [k(s)\sin\theta(s)]^\prime}{k(s)}
\right]^2 ds\ne 0.
\]
This ends the proof.
\end{proof}

\begin{remark}
When the curve $\mathbf{c}$ is such that $\int_0^L \tau(s)ds=2m\pi, \; m\in \mathbb Z\setminus \{0\}$
there are no ruled surfaces
as given by equation \eqref{eq:a}
containing $\mathbf{c}$ and having it as a closed principal
curvature line. In this situation we have always umbilic points
along $\mathbf{c}$.
In fact, in this case $k_2(s)= A(s) = 0$ and $m \neq 0$ implies that
$sin\theta(s) $ always vanishes. These points, at which $k_1(s)$ also
vanishes, happen to be the umbilic points. See Proposition \ref{prop:lc}.
\end{remark}

\begin{corollary}\label{cor:1}
Let $\mathbf{c}$ be a closed planar or spherical
Frenet curve of length $L$
in $\mathbb R^3$. Then there exists a germ of an oriented surface
containing $\mathbf{c}$ and having it as a hyperbolic principal
cycle if and only if $\mathbf{c}$ is not a circle.
\end{corollary}

\begin{proof}
In the case of a planar curve, let $\mathbf{c}(s)=(x(s),y(s),0)$ with
curvature $k$ and consider the Frenet frame $\{t,n,{\mathbf z}\}$,
${\mathbf z}=(0,0,1)$ associated to $\mathbf{c}$. Any parametrized surface
$\alpha$ containing $\mathbf{c}$ as a principal curvature line has the normal
vector equal to $N=\cos\theta(s) n(s) +\sin\theta(s) {\mathbf z}$.
Therefore,
\[
N^\prime = -k(s)\sin\theta(s)
t(s)+\theta^\prime[-\sin\theta(s){\mathbf z}+\cos\theta(s) n(s)].
\]

By Rodrigues formula $N^\prime=-\lambda(s) t(s)$ if and only if
$\theta(s)=\theta_0=cte$. One principal curvature is equal to
$k_1(s)=k(s)\sin\theta_0$. By the criterium of hyperbolicity of a
principal cycle, see Remark \ref{r:1} and Theorem \ref{th:lp}, the
principal curvatures can not be constant along a principal cycle. A
construction of the germ of surface containing $\mathbf{c}$ as a hyperbolic
principal cycle can be done as in Proposition \ref{prop:lc} and
Theorem \ref{th:lp2}. In the case of spherical curves, any closed
curve has total torsion equal to zero and the proof follows the same
steps of the planar case. This ends the proof.
\end{proof}

\section{Concluding Remarks}

The study of principal lines goes back to the works of Monge, see
\cite[page  95]{st}, Darboux \cite{da} and many others. In
particular, the local behavior of principal lines near umbilic
points is a classical subject of research, see \cite{gs3} for a
survey. The structural stability theory and dynamics of principal
curvature lines on surfaces was initiated by Gutierrez and Sotomayor
\cite{gs1}, \cite{gs2} and also  has been the subject of recent
investigation  \cite{arxiv}.

The possibility of a Frenet (biregular) closed curve in the space to
be a principal line of a surface along it depends only on its total
torsion to be an integer multiple of $2\pi.$

The presence of umbilic points on such surface depends on function
$A(s)$ as well as on $k(s)$ and $\theta$, which in turn depends on
$\tau$. In fact $\theta  = -\int _0 ^s \tau(s)ds + \theta_0$,
depending on a free parameter $\theta_0$.

In fact the location of the umbilic points given by
equation \eqref{eq:umb} change  with the parameter $\theta_0$.

In this paper it has been shown that a non circular closed Frenet curve
$\mathbf{c}$ in $\mathbb R^3$ can be a principal line of a germ of
surface provided its total torsion is $2m\pi,\; m\in\mathbb Z$. In
Theorem \ref{th:lp2} the germ of the surface has been constructed in such a way that
$\mathbf{c}$ is a hyperbolic principal line. As it is well known
that the total torsion of a closed curve can be any real number, the
results of this paper  show  that closed principal lines (principal
cycles) are special curves of $\mathbb R^3$.
This completes the results established in  \cite{pc} and \cite{li}.

The generic behavior of principal curvature  lines near a regular curve of umbilics was studied in \cite{gsabc}.

\vskip .7cm
{\bf Acknowledgements:}  This paper was done under the project CNPq 473747/2006-5.
  The first author had the  support of FUNAPE/UFG.

\vskip 2cm
\author{
\noindent Ronaldo A. Garcia,\\Instituto de Matem\'{a}tica e Estat\'{\i}stica,\\
Universidade Federal de Goi\'as,\\CEP 74001--970, Caixa Postal 131,
\\Goi\^ania, GO, Brazil.\\
E-mail: ragarcia@mat.ufg.br\\
\\
Luis F. Mello,\\
Instituto de Ci\^encias Exatas,\\
Universidade Federal de Itajub\'a,\\
CEP 37500--903, Itajub\'a, MG, Brazil.\\
E-mail: lfmelo@unifei.edu.br \\
\\
 Jorge Sotomayor\\Instituto de Matem\'{a}tica e Estat\'{\i}stica,\\Universidade de S\~{a}o Paulo,
\\Rua do Mat\~{a}o 1010, Cidade Universit\'{a}ria, \\CEP 05508-090, S\~{a}o Paulo, S.P., Brazil\\
E-mail: sotp@ime.usp.br}

\end{document}